\documentclass[draft]{amsart}
\usepackage[all,arc,curve,color,frame]{xypic}
\usepackage{amsfonts,amssymb,xy}
\theoremstyle{plain}
\newtheorem{theorem*}{Theorem}
\newtheorem*{lemma*}{Lemma}
\newtheorem{corollary*}{Corollary}
\newtheorem*{proposition*}{Proposition}
\newtheorem{conjecture*}{Conjecture}
\newtheorem{theorem}{Theorem}[section]
\newtheorem{lemma}[theorem]{Lemma}
\newtheorem{corollary}[theorem]{Corollary}
\newtheorem{proposition}[theorem]{Proposition}

\theoremstyle{remark}

\newtheorem*{remarks}{Remarks}

\newtheorem*{example}{Example}
\newtheorem*{claim}{Claim}

\theoremstyle{definition}

  \def\F{\Bbb{F}} \def\Z{\Bbb{Z}}  
\def\N{\Bbb{N}}   \def\ll{\langle} \def\rr{\rangle}
 \def\a{\alpha} \def\g{\gamma}  \def\bp{\begin{pmatrix}}
 \def\ep{\end{pmatrix}} \def\bn{\begin{enumerate}} 
   \def\en{\end{enumerate}}
\def\ba{\begin{array}} \def\ea{\end{array}}  
   \def\a{\alpha} \def\b{\beta} 
\def\id{\operatorname{id}} \def\Aut{\operatorname{Aut}} \def\im{\operatorname{Im}} 
  
\def\ker{\operatorname{Ker}}\def\be{\begin{equation}} \def\ee{\end{equation}}

     \def\fr12{\frac{1}{2}} \def\z12{\Z[\fr12]}

\def\ol{\overline}

\def\G{\Gamma}

\def\shtimes{\hspace{-0.1cm}\times\hspace{-0.1cm}}

\newcommand{\abs}[1]{\lvert#1\rvert}

\DeclareMathAlphabet{\mathbf}{OML}{cmm}{b}{it}

\numberwithin{equation}{section}

\begin{document}

\title[HNN Extensions of $p$-Groups which are Residually $p$ ]{A Criterion for HNN Extensions of Finite $p$-Groups to be Residually $p$}
\author{Matthias Aschenbrenner}
\address{University of California, Los Angeles, California, USA}
\email{matthias@math.ucla.edu}
\thanks{The first author was partially supported by a grant from the National Science Foundation.}

\author{Stefan Friedl}
\address{University of Warwick, Coventry, UK}
\email{s.k.friedl@warwick.ac.uk}

\date{June 7, 2010}
\begin{abstract}
We give a criterion for an HNN extension of a finite $p$-group to be residually $p$.
\end{abstract}

\maketitle

\section{Statement of the Main Results}\label{sec:Introduction}

\noindent
By an \emph{HNN pair} we mean a pair $(G,\varphi)$ where $G$ is a group and $\varphi\colon A\to B$ is an isomorphism between subgroups $A$ and $B$ of $G$. Given such an HNN pair $(G,\varphi)$ we consider the corresponding HNN extension
\[ G^* = \ll G,t \, |\, t^{-1}at = \varphi(a),\ a\in A\rr \]
of $G$, which we denote, by slight abuse of notation, as $G^* = \ll G,t \, |\, t^{-1}At=\varphi(A)\rr$.
Throughout this paper we fix a prime number $p$, and by a  $p$-group we mean a finite group of $p$-power order. We are interested in the question under which conditions an HNN extension of a $p$-group is residually a $p$-group. (HNN extensions of finite groups are always residually finite \cite{BT78, Co77}.)
Recall that given a property $\mathcal P$ of groups, a group $G$ is said to be \emph{residually $\mathcal P$} if for any non-trivial $g\in G$ there exists a morphism $\a\colon G\to P$ to a group $P$ with property $\mathcal P$ such that $\a(g)$ is non-trivial.

Given HNN pairs $(G,\varphi)$ and $(G',\varphi')$, a group morphism $\alpha\colon G\to G'$ is a morphism of HNN pairs if $\alpha(A)\subseteq A'$, $\alpha(B)\subseteq B'$, and the diagram
\[ \xymatrix{ A' \ar[r]^{\varphi'} & B' \\ A\ar[r]^{\varphi}\ar[u]^{\alpha|_A} & B \ar[u]_{\alpha|_B}}\]
commutes. (When talking about an HNN pair $(G,\varphi$),  we always denote the domain and codomain of $\varphi$ by $A$ respectively $B$, possibly with decorations.) A morphism $\alpha\colon (G,\varphi)\to (G',\varphi')$ of HNN pairs is  called an \emph{embedding of HNN pairs} if $\alpha$ is injective. Given a group $G$ and $g\in G$ we denote the conjugation automorphism $x\mapsto g^{-1}xg$ of $G$ by $c_g$.

There is a  well-known criterion for HNN extensions to be residually $p$:

\begin{lemma}\label{lem:resp}
Let $(G,\varphi)$ be an HNN pair, where $G$ is a $p$-group.
Then the following are equivalent:
\bn
\item the group $G^*=\ll G,t \, |\, t^{-1}At=\varphi(A)\rr$ is residually $p$;
\item there exists a $p$-group $X$ and an automorphism $\g$ of $X$ of $p$-power order such that $(G,\varphi)$ embeds into $(X,\gamma)$;
\item there exists a $p$-group $Y$ and $y\in Y$ such that $(G,\varphi)$ embeds into $(Y,c_y)$.
\en
\end{lemma}

\begin{proof}
For a proof of the equivalence of (1) and (3) we refer to \cite[Proposition~1]{RV91}.
Clearly (3) implies (2). On the other hand, let $X$ be a $p$-group and $\g\in\Aut(X)$ such that $(G,\varphi)$ embeds into $(X,\gamma)$, and suppose $\g$ has order $p^k$. Let
$Y=\Z/p^k\Z\ltimes X$ where $1\in \Z/p^k\Z$ acts on $X$ on the right via $\g$, and let $y=(1,1)\in \Z/p^k\Z \ltimes X$. Then $(G,\varphi)$ embeds into $(Y,c_y)$.
\end{proof}

\begin{example}
Suppose $A=B=G$, i.e., $G^*=\ll t\rr \ltimes G$ where $t$ acts on $G$ on the right via $\varphi$. Then $G^*$ is residually $p$ if and only if the automorphism $\varphi$ of $G$ has order a power of $p$.
\end{example}

Let $G$ be a group.
We say that a finite sequence  $\mathbf G = (G_1,\dots,G_n)$ of normal subgroups of $G$ with
$$G=G_1\supsetneq G_2\supsetneq\cdots\supsetneq G_n=\{1\}$$
is a \emph{filtration} of $G$. Given such a filtration $\mathbf G$ of $G$ we set $G_i:=\{1\}$ for $i>n$. We say that $\mathbf G$ is \emph{central} if $G_{i}/G_{i+1}$ is central in $G/G_{i+1}$ for each $i$.
Recall that the lower central series of $G$ is the sequence $(\gamma_i(G))_{i\geq 1}$ defined inductively by $\gamma_1(G)=G$ and $\gamma_{i+1}(G)=[G,\gamma_i(G)]$ for $i\geq 1$.
By definition $G$ is nilpotent if and only if $\gamma_n(G)=\{1\}$ for some $n\geq 1$. In that
case, taking $n$ minimal such that $\gamma_n(G)=\{1\}$, we obtain a central filtration  $\gamma(G)=(\gamma_1(G), \dots,\gamma_n(G))$ of $G$. In fact, $G$ admits a central filtration if and only if $G$ is nilpotent.

A filtration of a group $G$ is called a \emph{chief filtration} of $G$ if the filtration cannot be refined non-trivially to another filtration of $G$. It is well-known that a filtration $(G_i)$ of a $p$-group is a chief filtration if and only if all its non-trivial factors $G_i/G_{i+1}$ are of order $p$.
Note that a chief filtration of a $p$-group is necessarily a central filtration, since $\Z/p\Z$ has no non-trivial automorphism of $p$-power order.

We say that an HNN pair $(G,\varphi)$ and a filtration $(G_i)$ of $G$ as above are \emph{compatible} if $\varphi$ restricts to an isomorphism $A\cap G_i\to B\cap G_i$, for each $i$.
We recall the following theorem, which gives an intrinsic criterion for an HNN extension of a $p$-group to be residually $p$. This theorem
was shown in \cite[Lemma~1.2]{Ch94} (and later rediscovered in \cite{Mo07}); it
 can be viewed as a version for HNN extensions of Higman's theorem \cite{Hi64}, which gives a criterion for an amalgamated product of two $p$-groups to be  residually $p$.

\begin{theorem} \label{thm:hnnchief}
Let $(G,\varphi)$ be an HNN pair, where $G$ is a $p$-group.
Then the following are equivalent:
\bn
\item the group $G^*=\ll G,t \, |\, t^{-1}At=\varphi(A)\rr$ is residually $p$;
\item there exists a chief filtration $(G_1,\dots,G_n)$ of $G$, compatible with $(G,\varphi)$, such that
$\varphi(a) \equiv a \bmod G_{i+1}$ for all $i$ and $a\in A\cap G_i$.
\en
\end{theorem}

The objective of this paper is to give an alternative criterion for $G^*$ to be residually $p$, employing a certain group $H(G,\varphi)$ associated to every HNN pair $(G,\varphi)$, and defined as follows: set
$H_0=A\cap B$ and inductively $H_{i+1}=\varphi^{-1}(H_i)\cap H_i\cap \varphi(H_i)$ and put
$H(G,\varphi):=\bigcap_i H_i$.
Note that $\varphi$ restricts to an automorphism of $H(G,\varphi)$; in fact $H(G,\varphi)$ is the largest subgroup of $G$ which gets mapped onto itself by $\varphi$.
The group $H(G,\varphi)$ was introduced by Raptis and Varsos in \cite{RV89,RV91}. It had been previously employed in \cite{Hi81}, and a slight variant (the largest normal subgroup of $G$ mapped onto itself by $\varphi$) also occurs in \cite[Section~1.26]{Ba93}.
We propose to call $H(G,\varphi)$ the {\it core}\/ of $(G,\varphi)$. Indeed,  $H(G,\varphi)$ is the core  with respect to $\ll t\rr$ of $G$ construed as a subgroup of $G^*$, i.e., $H(G,\varphi)=\bigcap_{i\in\Z} t^{-i} G t^{i}$; so if $G$ is abelian, then $H(G,\varphi)$ is indeed the core of $G$ in $G^*$ (the largest normal subgroup of $G^*$ contained in $G$). See Lemma~\ref{lem:altdefh}, where we give other descriptions of the core of $(G,\varphi)$ which are oftentimes useful. Note that if $\alpha\colon (G,\varphi)\to (G',\varphi')$ is an embedding of HNN pairs, then $\alpha(H(G,\varphi)) \leq H(G',\varphi')$.

If $(G_i)$ is a filtration compatible with the  HNN pair $(G,\varphi)$, then  for any $i<j$ the morphism
$\varphi$ induces an isomorphism
$$(A\cap G_i)G_j/G_j \to (B\cap G_i)G_j/G_j,$$
which we denote by $\varphi_{ij}$. For $a\in A\cap G_i$, the conjugation automorphism $c_a$ of $G$ induces an automorphism of $(A\cap G_i)G_j/G_j$ which we continue to denote by $c_a$, and similarly for $c_b$ with $b\in B\cap G_i$. We can now formulate our first result, which gives an obstruction
to an HNN extension of a $p$-group being residually $p$. The statement of the proposition is inspired by the ideas of \cite{RV91}. {\it For the rest of this section we fix an HNN pair $(G,\varphi)$ where $G$ is a $p$-group, and we let $G^*=\ll G,t \, |\, t^{-1}At=\varphi(A)\rr$.}\/

\begin{proposition}\label{thm:obstruction}
If  $G^*$ is residually $p$, then
there exists a central filtration $\mathbf G=(G_i)$ of $G$, compatible with $(G,\varphi)$,
such that for any $i<j$, any $a\in A\cap G_i$ and any $b\in B\cap G_i$ the order of the automorphism  of  $H(G_i/G_j,c_b\circ \varphi_{ij} \circ c_a)$ induced by $c_b\circ \varphi_{ij} \circ c_a$ is a power of $p$.
\end{proposition}

In \cite[Theorem~13]{RV91} it is claimed that the following strong converse to Proposition~\ref{thm:obstruction} holds:
{\it If there exists a central filtration of $G$, compatible with $(G,\varphi)$, and if the order of the automorphism of $H(G,\varphi)$ induced by $\varphi$ is a power of $p$, then the HNN extension $G^*$ of $G$ is residually $p$.}\/
In Section~\ref{section:example} we show that this statement is incorrect; in fact, we give two counterexamples to
\cite[Theorem~13]{RV91}, highlighting the role of $a$ and $b$ and the importance of the filtration $\mathbf G$ in Proposition~\ref{thm:obstruction}.

\medskip

Our main theorem is  the following converse to Proposition~\ref{thm:obstruction}.

\begin{theorem}\label{thm:resp}
Suppose there exists a central filtration $\mathbf G=(G_i)$ of $G$, compatible with $(G,\varphi)$,
such that for any $i$ the order of the automorphism of  $H(G_{i}/G_{i+1},\varphi_{i,i+1})$ induced by $\varphi_{i,i+1}$
is a power of $p$. Then $G^*$ is residually $p$.
\end{theorem}

Note that for a chief filtration $\mathbf G$ of $G$, the statement of Theorem~\ref{thm:resp} is equivalent to the implication (2)~$\Rightarrow$~(1) in Theorem~\ref{thm:hnnchief}.

\medskip

For every filtration $(G_i)$ of $G$ compatible with $(G,\varphi)$ and any $i$, the group $H(G_{i}/G_{i+1},\varphi_{i,i+1})$
is a subgroup of $H(G/G_{i+1},\varphi_{i+1})$; here $\varphi_i$ is the isomorphism $AG_i/G_i\to BG_i/G_i$ induced by $\varphi$.
We therefore get the following corollary.

\begin{corollary}\label{cor:resp}
Assume that there exists a central filtration $(G_i)$ of $G$, compatible with $(G,\varphi)$,
such that for any $i$ the order of the automorphism  of  $H(G/G_{i},\varphi_i)$ induced by $\varphi_i$
is a power of $p$. Then $G^*$ is residually $p$.
\end{corollary}



\subsection*{Conventions} All groups are finitely generated. By a $p$-group we mean a finite group of $p$-power order. The identity element of a multiplicatively written group is denoted by $1$.

\subsection*{Acknowledgment} We would like to thank the anonymous referee for many helpful suggestions.

\section{Preliminaries on the Core of an HNN Pair} \label{section:h}

\noindent
In this section, we let $(G,\varphi)$ be an HNN pair with corresponding HNN extension $G^*$ of $G$. For $g\in G$ and $n\in\N$ we say that $\varphi^n(g)$ is defined if $g$ is in the domain of the $n$-fold compositional iterate of $\varphi$ thought of as a partially defined map $G\to G$, and similarly we say that $\varphi^{-n}(g)$ is defined if $g$ is in the domain of the $n$-fold iterate of $\varphi^{-1}$.
We first prove:

\begin{lemma} \label{lem:altdefh}
The group $H=H(G,\varphi)$ is the largest subgroup of $G$ such that $\varphi(H)=H$, and as subgroups of $G^*$,
\begin{equation}\label{equ:hdef0}
H = \bigcap_{i\in\Z} t^{-i} G t^{i}.
\end{equation}
Moreover,
\be \label{equ:hdef} H =\{ g\in G : \text{$\varphi^j(g)$ is defined for all $j\in \Z$}\}.\ee
If $A$ is finite, then there exists an integer $r\geq 0$ such that for any $s\geq r$ we have
\be \label{equ:hdef2} H =\{ g\in G : \text{$\varphi^j(g)$ is defined for $j=0,\dots,s$}\}.\ee
\end{lemma}

\begin{proof}
Recall that we introduced $H=\bigcap_i H_i$ as the intersection of the inductively defined descending sequence of subgroups
$$A\cap B=H_0\supseteq H_1\supseteq\cdots \supseteq H_i\supseteq H_{i+1}=\varphi^{-1}(H_i)\cap H_i\cap \varphi(H_i)\supseteq \cdots$$ of $G$. Clearly $\varphi(H_{i+1})\subseteq \varphi(\varphi^{-1}(H_i))=H_i$ for each $i$, hence $\varphi(H)\subseteq \varphi(\bigcap_i H_{i+1})\subseteq H$; similarly,
$\varphi^{-1}(H)\subseteq H$, hence $\varphi(H)=H$. Moreover,  given any $H'\leq G$ with $\varphi(H')=H'$, an easy induction on $i$ shows that $H'\subseteq H_i$ for all $i$, so $H'\subseteq H$.
To prove \eqref{equ:hdef0} and \eqref{equ:hdef} we show, by induction on $i$:
\be\label{equ:hdef3} H_{i} =  \bigcap_{\abs{j}\leq i+1} t^{-j} G t^{j}
    =  \{ g\in G : \text{$\varphi^j(g)$ is defined for all $j$ with $\abs{j}\leq i+1$}\}.
\ee
 For $i=0$ note that the Normal Form Theorem for HNN extensions \cite[Theorem~11.83]{Ro95} yields $G\cap tGt^{-1}=A=tBt^{-1}$, hence
$$t^{-1} G t \cap G \cap t G t^{-1} = A\cap B=H_0.$$
Moreover, given $g\in G$, clearly both $\varphi(g)$ and $\varphi^{-1}(g)$ are defined precisely if $g\in A\cap B$.
Now suppose \eqref{equ:hdef3} has been shown for some $i\geq 0$. Since $H_i\subseteq A\cap B$ we have
$$H_{i+1} = \varphi(H_i)\cap H_i\cap \varphi^{-1}(H_i) = t^{-1} H_i t\cap H_i\cap t H_i t^{-1} = \bigcap_{\abs{j}\leq i+1} t^{-j} G t^{j},$$
where we used the inductive hypothesis for the last equality. Now let $g\in G$. Then $g\in H_{i+1}$ if and only if $\varphi(g),g,\varphi^{-1}(g) \in H_i$. By the inductive hypothesis,  this in turn is equivalent to $\varphi^j(g)$ being defined for all $j$ with $\abs{j}\leq i+1$.

For $i\geq 0$ we now define
$$H_i' =\{ g\in G : \text{$\varphi^j(g)$ is defined for $j=0,\dots,i$}\}$$
and set $H':=\bigcap_i H_i'$. Then clearly $H\subseteq H'$ and $\varphi(H')\subseteq H'$. Suppose $A$ is finite; then there is an integer $r\geq 0$ with $H'=H_r'=H_{r+1}'=\cdots$. To show \eqref{equ:hdef2} we prove that $\varphi(H')=H'$ (which yields $H=H'$ by the first part of the lemma). Let $g\in H'$. Since $A$ is finite, there exist $k\geq 0$ and $l>0$ such that $\varphi^{k+l}(g)=\varphi^k(g)$. Since $\varphi$ is injective, $\varphi^l(g)=g$, hence $\varphi^{-1}(g)=\varphi^{l-1}(g)$ exists, and clearly $\varphi^{-1}(g)\in H'$. Hence $\varphi^{-1}(H')\subseteq H'$ and thus $\varphi(H')=H'$ as required.
\end{proof}

The main difficulty in dealing with the core of $(G,\varphi)$ is that it does not behave well under taking quotients.
For example,  if  $H(G,\varphi)$ is trivial, and if $K\leq G$
such that $\varphi(A\cap K)=B\cap K$ and $\ol{\varphi}\colon AK/K\to BK/K$ is the isomorphism induced by $\varphi$, then it is not true in general that $H(G/K,\ol{\varphi})$ is trivial.
A non-abelian example of this phenomenon is given in Section~\ref{section:ex2}.
But this can even happen in the abelian case:

\begin{example}
Suppose
$$G=\F_p^3,\quad  A=\big\{ (a_1,a_2,0): a_i\in \F_p\big\}, \quad B=\big\{ (b_1,0,b_2): b_i\in \F_p\big\}.$$
Then $A\cap B=\big\{(c,0,0): c\in \F_p\big\}$. Let $x,y,z\in \F_p$ with $x,z\ne 0$, and
suppose $\varphi$ is the isomorphism $A\to B$ given by
$$\varphi(a_1,a_2,0)=(xa_1,0,ya_1+za_2).$$
If $y\ne 0$, then $\varphi(A\cap B)\cap (A\cap B)=0$, in particular $H(G,\varphi)=0$.
Now let $$K=\big\{(0,k_1,k_2) : k_1,k_2\in \F_p\big\}.$$
Then $\varphi(A\cap K)=B\cap K$,
but $H(G/K,\ol{\varphi})$ is non-trivial, in fact it equals $G/K\cong \F_p$ and the automorphism induced by $\varphi$ is multiplication by $x$.
\end{example}

A very similar example shows that the proof of \cite[Lemma~6]{RV91} does not work in general:

\begin{example}
Suppose
$$G=\F_p^4,\quad  A=\{ (a_1,a_2,a_3,0): a_i\in \F_p\}, \quad B=\{ (b_1,0,b_2,b_3): b_i\in \F_p\},$$
so $A\cap B=\{(c_1,0,c_2,0):c_i\in\F_p\}$. Let $a,b,c\in \F_p$ with $a\neq 0$, and suppose $\varphi$ is the isomorphism $A\to B$ given by
$$\varphi(a_1,a_2,a_3,0)=(aa_1,0,ba_1+a_2,ca_1+a_3).$$
If $c\neq 0$, then $H=0$, and for $x=(0,0,1,0)\in A\cap B$ we have $L_x=\{(0,d_1,d_2,d_3):d_i\in\F_p\}$ (employing the notation of \cite{RV91}), but $H(G/L_x,\ol{\varphi})\neq 0$, contrary to what is assumed in the inductive step in the proof of  \cite[Lemma~6]{RV91}.
\end{example}

\section{Obstructions to HNN Extensions Being Residually $p$}
\label{section:firstproofs}

\noindent
Before we give the proof of Proposition~\ref{thm:obstruction} we prove the following lemma. Again, we let $(G,\varphi)$ be an HNN pair.

\begin{lemma} \label{lem:powerofp}
Let $Y$ be a group, $y\in Y$, and  $\a\colon (G,\varphi)\to (Y,c_y)$ an embedding of HNN pairs.
\bn
\item If $Y$ is a $p$-group, then the order of the restriction of $\varphi$ to $H(G,\varphi)$ is a power of $p$.
\item Let  $a\in A$, $b\in B$; then $\a$ is an
embedding  $(G,c_b\circ \varphi\circ c_a)\to (Y,c_{ayb})$.
\en
\end{lemma}

\begin{proof}
For the first statement write $H=H(G,\varphi)$ and identify $G$ with $\a(G)\leq Y$.
Then $\varphi|_H=c_y|_H$, hence the order of $\varphi|_H$ divides the order of $y$.
The second statement follows immediately from $c_{ayb}=c_b\circ c_y\circ c_a$.
\end{proof}

\begin{proof}[Proof of Proposition~\ref{thm:obstruction}]
Suppose $G$ is a $p$-group and $\ll G,t \, |\, t^{-1}At=\varphi(A)\rr$ is residually $p$.
By Lemma~\ref{lem:resp} we can find a $p$-group $Y$, $y\in Y$ and an embedding $\alpha\colon (G,\varphi)\to (Y,c_y)$  of HNN pairs. We identify $G$ with its image under $\alpha$.

Let $(Y_i)$ be any central filtration of $Y$, and let
$G_i=Y_i\cap G$ for each $i$.
Evidently $G_{i}/G_{i+1}$ is central in $G/G_{i+1}$ for any $i$.
Possibly after renaming we can also achieve that for each $i$ we have $G_{i+1}\subsetneq G_{i}$,
i.e., $\mathbf G=(G_i)$ is a central filtration of $G$. Furthermore note that for any $i$ the following holds:
\[ \varphi(A\cap G_i)=\varphi(A\cap G\cap Y_i)=c_y(A\cap Y_i)=B\cap Y_i=B\cap G_i\]
since $Y_i$ is normal in $Y$. This shows that $\mathbf G$ is compatible with $(G,\varphi)$.

Finally let $i<j$. Then  $\a$ gives rise to  an embedding
$$(G_i/G_j,\varphi_{ij})\to (G/G_j,\varphi_j)\to (Y/Y_j,c_{yY_j})$$
of HNN pairs. It follows now from Lemma~\ref{lem:powerofp} that for any $a\in A\cap G_i$, $b\in B\cap G_i$ the order of the restriction of  $c_b\circ \varphi_{ij}\circ c_a$ to $H(G_i/G_j,\varphi_{ij})$
is a power of $p$.
\end{proof}

\section{Examples} \label{section:example}

\noindent
In this section we apply Proposition~\ref{thm:obstruction} to two HNN extensions.
The first example highlights the role of $a$ and $b$ in Proposition~\ref{thm:obstruction}, the second one shows the importance of the central series.
Both are counterexamples to \cite[Theorem~13]{RV91}.

\subsection{The first example} \label{section:ex1}

The multiplicative group $$P:=\ll x,y\, |\, x^3=y^3=[x,y]=e\rr$$
is naturally isomorphic to the additive group $\F_3\oplus \F_3$.
We write $\ll x \rr$ and $\ll y \rr$ for the subgroups of $P$ generated by $x$ and $y$, respectively.
We think of elements in the group ring $\F_3[P]$ as polynomials $f(x,y)=\sum_{i,j=0}^2 v_{ij}x^iy^j$ with coefficients $v_{ij}\in\F_3$.
Furthermore $f(x)$ always denotes an element in the subring $\F_3[\ll x\rr]$ of $\F_3[P]$ and similarly $f(y)$ will denote an element in $\F_3[\ll y \rr]\subseteq \F_3[P]$.

Let $G=P\ltimes \F_3[P]$ where $P$ acts on its group
ring $\F_3[P]$ by multiplication. (Here $P$ is a multiplicative group and $\F_3[P]$ is an additive group. Note that $G$ is in fact just
the wreath product $\F_3 \wr P$.) Evidently $G$ is a $3$-group.
For $f\in \F_3[P]$ we have
\[ c_{(x^ny^m,0)}(1,f(x,y))=(x^{-n}y^{-m},0)(1,f(x,y))(x^{n}y^{m},0)=(1,x^{n}y^{m}f(x,y)).\]
Now consider the subgroups
$A=\ll x \rr \ltimes \F_3[\ll x \rr]$ and $B=\ll y \rr \ltimes \F_3[\ll y \rr]$ of $G$, and let $\varphi\colon A\to B$ be the map given by
\[ \varphi(x^n,f(x))=(y^n,2y^{-1}f(y)).\]
It is straightforward to verify that $\varphi$ is indeed an isomorphism. In fact,  $\varphi$ is the restriction to $A$ of the automorphism $\phi$ of $G=P\ltimes \F_3[P]$ given by
$$(x^ny^m,f(x,y))\mapsto  (x^my^n,2y^{-1}f(y,x)).$$

\begin{claim}
The HNN pair $(G,\varphi)$ is compatible with the lower central series $\g(G)$ of $G$.
\end{claim}

The claim follows immediately from the fact that $\varphi$ is the restriction of an automorphism of $G$, and the fact that the groups in the lower central series are characteristic. Indeed, we compute
\[ \ba{rcl} \varphi(A\cap \gamma_i(G))&=&\phi(A\cap \gamma_i(G))=\\
&=&\phi(A)\cap \phi(\gamma_i(G))=\phi(A)\cap \gamma_i(G)=B\cap \gamma_i(G).\ea \]

\begin{claim}
The subgroup $H(G,\varphi)$ of $G$ is trivial.
\end{claim}

Indeed, first note that $A\cap B=\{ (1,v) : v\in \F_3\}$.
But $\varphi(1,v)=(1,2y^{-1}v)$. This shows that $(A\cap B)\cap \varphi(A\cap B)=\{1\}$, hence $H(G,\varphi)=\{1\}$.

\medskip

If \cite[Theorem~13]{RV91} was correct, then $\ll G,t \, |\, t^{-1}At= \varphi(A)\rr$ would have to be a group which is residually a $3$-group.
But the combination of the next claim with Proposition~\ref{thm:obstruction} shows that this is not the case.

\begin{claim}
Put $\psi:=\varphi \circ c_{(x,0)}\colon A\to B$.
Then $H(G,\psi)\neq\{1\}$,
and the restriction of $\psi$ to  $H(G,\psi)$  has order $2$.
\end{claim}

We have
$$\psi(1,v)=(\varphi \circ c_{(x,0)})(1,v)=\varphi(1,vx)=(1,2y^{-1}vy)=(1,2v).$$
This shows that $\psi$ induces an automorphism of $A\cap B$, and the automorphism has  order $2$. It follows immediately
 that $H(G,\psi)=\{ (1,v) : v\in \F_3\}$ and that $\psi$ restricted to  $H(G,\psi)$ has order $2$.

\subsection{The second example}\label{section:ex2}
In the following we write elements of $\F_3 \oplus \F_3 \oplus \F_3$ as column vectors. The automorphism of $\F_3 \oplus \F_3 \oplus \F_3$ given by the matrix
\[ X:=\bp 0&1&0 \\ 0&0&1 \\ 1&0&0 \ep\]
clearly descends to an automorphism of
$$V := (\F_3 \oplus \F_3 \oplus \F_3 ) / \big\{(a,a,a)^{\operatorname{t}}:a\in \F_3\big\}.$$
In the rest of this section a column vector in $\F_3\oplus \F_3 \oplus \F_3$ will always stand for the element in $V$ it represents. Consider the $3$-group
$G:=\ll x\, |\, x^3\rr \ltimes V$ where $x$ acts on $V$ on the right via $X$.
Note that given integers $m$, $n$ and $u,v\in V$ we have
\[ \ba{rcl} [(x^m,u),(x^n,v)]&=&(x^m,u)\cdot(x^n,v)\cdot(x^{-m},-X^{-m} u)\cdot(x^{-n},-X^{-n} v)\\ &=& \big(1,X^{-n}(X^{-m} v-v)-X^{-m}(X^{-n} u-u)\big).\ea \]
Since for $r=1,2$ we have
$$(X^r-\id)(V)=\big\{(w_1,w_2,w_3)^{\operatorname{t}} : w_1+w_2+w_3=0\big\},$$
it follows that
\[ \ba{rcl} \gamma_2(G)=[G,G]&=&\big\{ (1,w) : w\in (X-\id)(V) \big\}\\
&=&\big\{ \big(1,(w_1,w_2,w_3)^{\operatorname{t}}\big) : w_1+w_2+w_3=0\big\}.\ea \]
A similar calculation shows that
\[ \gamma_3(G)=[G,[G,G]]=\big\{ \big(1,(a,a,a)^{\operatorname{t}}\big) : a\in \F_3\big\}=0.\]
Now let
\[ a=\bp 1\\ 0 \\ 0\ep \text{ and } b=\bp 1\\1\\-1 \ep \in V\subseteq G.\]
Note that $a$ and $b$ represent the same element in $G/[G,G]$.
Let $A$ and $B$ be the subgroups of $V\leq G$ generated by $a$ and $b$, respectively.
Let $\varphi\colon A\to B$ be the isomorphism given by $\varphi(a)=2b$.
Note that $A\cap \gamma_2(G)=B\cap \gamma_2(G)=\{1\}$.
It follows that the HNN pair $(G,\varphi)$ is compatible with the filtration of $G$ given by the lower central series.
Finally note that $A$ and $B$ intersect trivially. In particular $H(G,\varphi)$ is trivial.

If \cite[Theorem~13]{RV91} was correct, then $\ll G,t \, |\, t^{-1}At= \varphi(A)\rr$ would have to be a residually $3$-group.
The following lemma in conjunction with Proposition~\ref{thm:obstruction} shows that this is not the case.

\begin{lemma}
There exists no  central filtration $(G_i)$ of $G$, compatible with the HNN pair $(G,\varphi)$,
such that for any $i$ the order of  the automorphism of $H(G/G_i,\varphi_i)$ induced by $\varphi_i\colon AG_i/G_i \to BG_i/G_i$
is a power of $3$.
\end{lemma}

\begin{proof}
Let $G=G_1\supsetneq G_2\supsetneq G_3 \supsetneq\cdots \supsetneq G_n=\{1\}$ be a central filtration of $G$ compatible with $(G,\varphi)$.
Denote the natural surjection $G\to G/G_i$ by $\pi_i$.
Note that $\pi_2$ factors through $G/[G,G]$ and therefore $\pi_2(a)=\pi_2(b)$.

First assume that $\pi_2(a)\ne 0$. In that case the subgroup $\pi_2(A)$  of $G/G_2$ is isomorphic to $\F_3$, and $\pi_2(A)=\pi_2(B)$.
It follows that $H(G/G_2,\varphi_2)=\pi_2(A)=\pi_2(B)\cong \F_3$.
Furthermore, since $\pi_2(a)=\pi_2(b)$ and $\varphi(a)=2b$ it follows that the automorphism of
$H(G/G_2,\varphi_2)\cong \F_3$ induced by $\varphi_2$ is multiplication by $2$, which has order $2$, hence is not a power of $3$.

Now assume that $\pi_2(a)=0$.
In that case we have $G\supsetneq G_2\supseteq [G,G]$ and $a\in G_2$. Recall that $G/[G,G]\cong \F_3^2$ and $a\not\in [G,G]$.
It follows easily that $G_2=\{1\}\times V\subseteq G=\ll x\, |\, x^3\rr \ltimes V$.
Recall that by the definition of a central series, $G_2/G_3$ is central in $G/G_3$.
In particular we have $(1,Xv-v)=[(x^{-1},0),(1,v)]\in G_3$ for any $v\in V$. Also note that $G_2\supsetneq G_3$.
It now follows immediately  that
$$G_3=\big\{ \big(1,(w_1,w_2,w_3)^{\operatorname{t}}\big) : w_1+w_2+w_3=0\big\}=[G,G].$$
We have $\pi_3(a)=\pi_3(b)\ne 0$. We now apply the same argument as above to see
that $H(G/G_3,\varphi_3)\cong \F_3$ and that the order of the automorphism of
$H(G/G_3,\varphi_3)$ induced by $\varphi_3$ is not a power of $3$.
\end{proof}

\section{Proof of Theorem~\ref{thm:resp}}

\noindent
In subsection \ref{section:extn} we first establish a useful consequence of Theorem~\ref{thm:hnnchief}.
In subsection~\ref{sec:special} we then prove some special cases of Theorem~\ref{thm:resp}, and  we give the proof of the general case of this theorem in subsection~\ref{sec:general}.

\subsection{An extension lemma}\label{section:extn}
The following lemma
will play a prominent role in our proof of Theorem~\ref{thm:resp}.

\begin{lemma}\label{lem:extension}
Let $(G,\varphi)$ an HNN pair, where $G$ is a $p$-group.
Suppose there exists a  central filtration $\mathbf G=(G_i)$ compatible with $(G,\varphi)$ such that for each $i$ there is a $p$-group $Q_i$ containing  $L_i:=G_{i}/G_{i+1}$
\textup{(}the ``$i$-th layer'' of the filtration $\mathbf G$\textup{)} as a subgroup
such that the isomorphism
$$A_i:=(A\cap G_i)G_{i+1}/G_{i+1}\xrightarrow{\varphi_{i,i+1}} B_i:=(B\cap G_{i})G_{i+1}/G_{i+1}$$ between subgroups of $L_i$ induced by $\varphi$
is the restriction of an inner automorphism of $Q_i$. Then the HNN extension
$G^*=\ll G,t \, |\, t^{-1}At = \varphi(A)\rr$ of $G$ is residually $p$.
\end{lemma}

\begin{proof}
Throughout the proof we write $\varphi_i=\varphi_{i,i+1}$. (This differs from the use of this notation in Section~\ref{sec:Introduction}.)
Note that by Lemma~\ref{lem:resp} and Theorem~\ref{thm:hnnchief} we can take for each $i$
a chief filtration $(H_{i1},\dots,H_{im_i})$ of  $L_i$ such that $\varphi_i(A_i\cap H_{ij})=B_i\cap H_{ij}$ for any $j$
 and such that
\be \label{equ:varphia} {\varphi}_i(a)\equiv a\bmod H_{i,j+1}\qquad\text{ for all $j$ and $a\in A_i\cap H_{ij}$.} \ee
For each $i$ let $\pi_i\colon G_i\to G_i/G_{i+1}=L_i$ be the natural epimorphism.
We set $G_{ij}:=\pi_i^{-1}(H_{ij})\leq G_i$. For each $i$ we have $G_{i1}=G_i$ and $G_{in_i}=G_{i+1}$.

\begin{claim}
For any $i,j$ the subgroup $G_{ij}$ is normal in $G$.
\end{claim}

Denote by $\pi$ the natural surjection $G\to G/G_{i+1}$. Note that if we consider $L_i=G_i/G_{i+1}$ as a subgroup of $G/G_{i+1}$ as usual, then $\pi_i$ is the restriction of $\pi$ to $G_i$, so $G_{ij}=\pi^{-1}(H_{ij})$.
It therefore suffices to show that $H_{ij}$ is normal in $G/G_{i+1}$.
But this follows immediately from the fact that  $G_i/G_{i+1}$ lies in the center of $G/G_{i+1}$. This concludes the proof of the claim.

\medskip

We now get the following filtration of $G$:
\[G=G_{11}\supsetneq G_{12}\supsetneq\cdots \supsetneq G_{1m_1}=G_2=G_{21} \supsetneq G_{22}\supsetneq\cdots\supsetneq G_{nm_n}=\{1\}.\]
Evidently each successive non-trivial quotient is isomorphic to $\Z/p\Z$, hence the above is a chief filtration for $G$.
Finally note that \eqref{equ:varphia} implies that
$$\varphi(a)\equiv a \bmod G_{i,j+1}\qquad\text{ for all $i$,  $j$ and $a\in A\cap  G_{ij}$.}$$
Hence the chief filtration satisfies condition (2) in Theorem~\ref{thm:hnnchief}, and we conclude that  the HNN extension
$G^*$  is residually $p$.
\end{proof}

\subsection{Two special cases of Theorem~\ref{thm:resp}}\label{sec:special}
The following proposition provides a proof of Theorem~\ref{thm:resp} in the case where $G$ is an elementary abelian $p$-group equipped with the trivial filtration and where furthermore $H(G,\varphi)=A\cap B$:

\begin{proposition}\label{prop:hnnembed}
Let $(G,\varphi)$ be an HNN pair with $G$ an elementary abelian $p$-group.
Suppose that $\varphi(A\cap B)=A\cap B$, i.e., $\varphi$ induces an automorphism of $A\cap B$, and assume that the order of this automorphism
is a power of $p$. Then there exists an elementary abelian $p$-group $X$ and an automorphism $\gamma$ of $X$ of $p$-power order such that  $(G,\varphi)$ embeds into $(X,\gamma)$.
In particular  $G^* =\ll G,t \, |\, t^{-1}at = \varphi(a),\ a\in A\rr$ is residually $p$.
\end{proposition}

Our proof of this proposition is inspired by the proof of Lemma~5 in \cite{RV91}. Below we use additive notation for abelian $p$-groups.

\begin{proof}
Choose $P\leq A$ and $Q\leq B$ such that $A=(A\cap B)\oplus  P$ and $B=(A\cap B)\oplus Q$, and then $S\leq G$ such that
$$G = {A}\oplus Q\oplus S=(A\cap B)\oplus P\oplus Q\oplus S.$$
We now let
$$X = {A} \oplus Q \oplus Q_1 \oplus\cdots\oplus Q_{p-2}\oplus S,$$
where each $Q_i=Q\times i$ is an isomorphic copy of $Q$.
We view $G$ as a subgroup of $X$ in the obvious way.
Let $\gamma$ be the endomorphism of $X$ such that for all $x\in X$,
$$\gamma(x) =
\begin{cases}
\varphi(x)     &  \text{if $x\in{A}$,}  \\
x\times 1      &  \text{if $x\in Q$,} \\
y\times (i+1)  & 	 \text{if $x=y\times i$ where $y\in Q$, $i=1,\dots,p-3$,} \\
\varphi^{-1}(y)&  \text{if $x=y\times (p-2)$ where $y\in Q$,}\\
x              &  \text{if $x\in S$.}
\end{cases}$$
By construction  the restriction of $\gamma$ to $A=(A\cap B)\oplus P$ equals $\varphi$. It remains to show that $\gamma$ is an automorphism of $p$-power order. By assumption  $\varphi|_{{A\cap B}}$ is an automorphism of ${A\cap B}$ of $p$-power order. Moreover, $\gamma$ restricts to an automorphism of the subgroup
$P\oplus Q \oplus Q_1\oplus\cdots\oplus Q_{p-2}$ of $X$ having order $p$, and $\gamma|_S=\id$. Hence $\gamma$ is an automorphism of $X$ of $p$-power order.
The last statement is now an immediate consequence of Lemma~\ref{lem:resp}.
\end{proof}

We now generalize the above proposition to the case  that $G$ is any abelian $p$-group.
We will again show that the resulting HNN~extension is residually $p$, but our construction will be somewhat less explicit in this case.

\begin{proposition}\label{prop:hnnresp}
Let $(G,\varphi)$ be an HNN pair with $G$ an abelian $p$-group.
Suppose that  $\varphi$ restricts to an automorphism of $A\cap B$ having $p$-power order. Then   $G^* =\ll G,t \, |\, t^{-1}at = \varphi(a),\ a\in A\rr$ is residually $p$.
\end{proposition}

\begin{proof}
First take a group embedding $\iota\colon G\to H:=(\Z/p^k\Z)^d$, for some $d$ and $k>0$ (cf.~\cite[p.~172]{RV91}). Then $\iota$ is an embedding $(G,\varphi)\to (H,\psi)$ of HNN pairs, where $\psi:=\iota\circ \varphi\circ \iota^{-1}\colon \iota(A) \to \iota(B)$.
After passing from $(G,\varphi)$ to $(H,\psi)$, we can therefore assume that $G=(\Z/p^k\Z)^d$.

Given $i>0$ we now write $G_i:=p^{i-1}G$. Note that for any subgroup $H\leq G$  the group $H\cap G_i$ is the \emph{characteristic} subgroup
$$H\cap G_i=\{ h\in H : p^{k-i+1}h=0 \}$$
of $H$. We thus get an (obviously central) filtration
$$G=G_1\supsetneq G_2\supsetneq\cdots\supsetneq G_k \supsetneq G_{k+1}=\{0\}$$
of $G$. Note that $G_k\cong (\Z/p\Z)^d$ and that for each $i\in \{1,\dots,k\}$ the isomorphism $\varphi$ restricts to an isomorphism between the characteristic subgroups $A\cap G_i\leq A$ and $B\cap G_i\leq B$.
Put differently, the filtration is compatible with $(G,\varphi)$.

\begin{claim}
There exists a $p$-group $Y$ and $y\in Y$ such that the HNN pair $(G_k,\varphi|_{A\cap G_k})$ embeds into $(Y,c_y)$.
\end{claim}

We write $A_k=A\cap G_k$ and $B_k=B\cap G_k$. Then
\[ (A\cap B)\cap G_k=(A\cap G_k)\cap (B\cap G_k)=A_k\cap B_k,\]
and $ (A\cap B)\cap G_k$ is characteristic in $A\cap B$. Recall that we assumed  that $\varphi$ restricts to an automorphism of $A\cap B$ of $p$-power order. It now follows easily from the above discussion
that $\varphi$ restricts to an automorphism of $A_k\cap B_k$  of $p$-power order.
Since $G_k$ is an elementary abelian $p$-group, the claim now follows immediately from  Proposition
\ref{prop:hnnembed}.

\medskip

For $i\in \{1,\dots,k\}$ we now write  $L_i:=G_{i}/G_{i+1}$; then $\varphi$ induces an isomorphism
\[A_i:=\big((A\cap G_i)+G_{i+1}\big)/G_{i+1}\xrightarrow{\varphi_{i}} B_i:=\big((B\cap G_{i})+G_{i+1}\big)/G_{i+1}\]
between subgroups of $L_i$. Note also that $$g\mapsto p^{k-i}g\colon G\to G$$ induces an isomorphism  $\Phi_i\colon L_i\to G_{k}$, giving rise to the following commutative diagram
\[ \xymatrix{ A_i\ar[d]_{\Phi_i|_{A_i}} \ar[r]^{\varphi_i} & B_i \ar[d]^{\Phi_i|_{B_i}} \\
A_k\ar[r]^{\varphi_k} & B_k,}\]
where the two vertical maps are injective.
In other words, $\Phi_i$ is an embedding of HNN pairs
$(L_i,\varphi_i)\to (L_k,\varphi_k)$.
By the above claim there exists a $p$-group $Y$ and $y\in Y$ such that
$(L_k,\varphi_k)=(G_k,\varphi|_{A_k})$ embeds into $(Y,c_y)$.
Hence $(L_i,\varphi_i)$ embeds into $(Y,c_y)$, for $i=1,\dots,k$. It now follows from Lemma~\ref{lem:extension}
that  the HNN extension
$G^*$ of $G$ is residually $p$.
\end{proof}

\subsection{The conclusion of the proof of Theorem \ref{thm:resp}.}
\label{sec:general}


We first formulate another interesting special case of Theorem~\ref{thm:resp}:

\begin{theorem}\label{thm:respabelian}
Let  $(G,\varphi)$  be an HNN pair where $G$ is an abelian $p$-group. Then the HNN extension $G^*=\ll G,t \, |\, t^{-1}At=\varphi(A)\rr$ of $G$ is residually $p$ if and only if the order of the restriction of $\varphi$ to $H(G,\varphi)$
is a power of $p$.
\end{theorem}

\begin{remarks} \mbox{}
\bn
\item
The HNN extension $G^*$ of $G$ is $\Z$-linear ~\cite[Corollary~3.5]{MRV08} and hence has, for each prime $q$, a finite-index subgroup which is residually $q$. In particular, $G^*$ always has a finite-index subgroup which is residually $p$.
\item This theorem also appears as \cite[Theorem~8]{RV91}. The proof in \cite{RV91} relies on the erroneous Lemma~6 in \cite{RV91}, but the referee informs us that the proof of \cite[Theorem~8]{RV91}
can be fixed.
\en
\end{remarks}

Assuming this theorem for a moment, we are now ready to prove Theorem~\ref{thm:resp} in general:

\begin{proof}[Proof of Theorem \ref{thm:resp}]
Let $(G,\varphi)$ be an HNN pair where $G$ is a $p$-group.
Assume that there exists a central filtration $(G_i)$ of $G$ compatible with $(G,\varphi)$,
such that for any $i$ the order of the automorphism of  $H(G_i/G_{i+1},\varphi_{i,i+1})$ induced by $\varphi$ is a power of $p$.
By Theorem~\ref{thm:respabelian} and Lemma~\ref{lem:resp}, for each $i$ there is an extension of $\varphi_{i,i+1}$ to an inner automorphism of a $p$-group containing $G_{i}/G_{i+1}$ as a subgroup.
Now Lemma~\ref{lem:extension} yields that $G^*$ is residually $p$.
\end{proof}

For the forward direction in Theorem~\ref{thm:respabelian} note that if $G^*$ is residually $p$, then so is its subgroup $\ll t\rr\ltimes H(G,\varphi)$, hence the order of $\varphi|_{H(G,\varphi)}$ is a power of $p$ by the example following Lemma~\ref{lem:resp}.
The remainder of this section will be occupied by the proof of the backward direction in Theorem~\ref{thm:respabelian}.
Throughout this section let  $(G,\varphi)$  be an HNN pair such that $G$ is an abelian $p$-group. 
In light of Lemma~\ref{lem:altdefh} there exists an integer $r\geq 1$
such that for any $s>r$ and $g\in G$ the following holds:
\be \label{equ:props}  g\in H(G,\varphi) \, \, \Longleftrightarrow \,\, \text{$\varphi^i(g)$ is defined for $i=0,\dots,s-1$.} \ee
For the time being  let $s$ be any integer such that $s>r$.
Consider the morphism
$$\phi\colon G^*=\ll G,t \, |\, t^{-1}At= \varphi(A)\rr \to \Z/s\Z
\quad\text{with $\phi(t)=1$ and $\phi(g)=0$ for $g\in G$.}$$
For $i\in \Z/s\Z$ we now write
$$G_i:=G\times i,\quad A_i:=A\times i,\quad B_i:=B\times i,$$
and we consider the groups
 $$K:=G_0*_{A_0=B_1} G_1*_{A_1=B_2}G_2 \cdots \ast_{A_{s-2}=B_{s-1}} G_{s-1}$$
 and
 \[  \G:=\ll K,t \, |\, t^{-1}A_{s-1}t=B_{0}\rr.\]
 We now have the following well-known lemma.

 \begin{lemma}
 The map
 \[ \ba{rcl} \G&\to & G^* \\
 (g,i) &\mapsto & t^{i}gt^{-i} \mbox{ for $g\in G$ and $i=0,\dots,s-1$}, \\
 t&\mapsto &t^s \ea \]
 defines an isomorphism onto $\ker(\phi)\leq G^*$.
 \end{lemma}

 \begin{proof}
 Note that the group $\G$ can be viewed as the fundamental group of a graph of groups based on a circuit of length $s$.
 On the other hand $G^*$ can be viewed as the fundamental group of a graph of groups based on a circuit of length $1$.
 The claim now follows from the theory of graphs of groups (see for example \cite{Se80}).
 Alternatively the lemma can also be proved easily using covering space theory. We leave the details to the reader.
 \end{proof}

We now consider $G':=H_1(K;\Z)$.
By a Mayer-Vietoris argument (cf.~\cite[II.2.8]{Se80}) we obtain an exact sequence
\be \label{equ:mv}\ba{cccccccccccccc}\bigoplus\limits_{i=0}^{s-2}&A_i&
\xrightarrow{\beta} &\bigoplus\limits_{i=0}^{s-1} G_i&\to& H_1(K;\Z)&\to&0\\
&a\shtimes i&\mapsto & a\shtimes i-\varphi(a)\shtimes (i+1)&&&&\ea \ee
where the morphism in the middle extends the morphisms $G_i=H_1(G_i;\Z)\to H_1(K;\Z)$ induced by the natural inclusions $G_i\to K$.

\begin{lemma}\label{lem:giinj}
Let $j\in \Z/s\Z$. The morphism $G_j\to H_1(K;\Z)$ is injective.
\end{lemma}

\begin{proof}
Let $a\in \bigoplus_{i=0}^{s-2}A_i$.
Write  $a=\sum_{i=i_1}^{i_2}a_i\shtimes i$ where $0\leq i_1\leq i_2\leq s-2$ and $a_i\in A$ for $i=i_1,\dots,i_2$, with $a_{i_1},a_{i_2}\ne 0$.
It is evident that the projections of $\b(a)$ to $G_{i_1}$ respectively $G_{i_2+1}$ are both non-zero (since $s\geq 2$).
It follows that $\im\beta\cap G_j=0$.
Exactness of \eqref{equ:mv} now yields that $G_j\to H_1(K;\Z)$ is injective.
\end{proof}

Note that $G'=H_1(K;\Z)$ is the quotient of an abelian $p$-group, in particular $G'$ itself is an abelian $p$-group.
Since $G_j\to G'$ is injective for any $j$ we can view
$A':=A_{s-1}$ and $B':=B_0$
 as subgroups of $G'$.
We denote by $\varphi'$ the isomorphism $A'\to B'$ defined by $\varphi$, i.e., $\varphi'(a\shtimes (s-1))=\varphi(a)\shtimes 0$ for all $a\in A$.

\medskip

The following lemma gives in particular a reinterpretation of $H(G,\varphi)$.

\begin{lemma}\label{lem:abprime}
\mbox{}
\bn
\item $H(G',\varphi')=A'\cap B'$.
\item Let $\Psi\colon G\xrightarrow{\cong} G_0\leq G'$ be the isomorphism given by $\Psi(g)=g\shtimes 0$ for $g\in G$. Then
$\Psi$ restricts to an isomorphism
\[H(G,\varphi)\to H(G',\varphi'),\]
and the following diagram commutes:
\[ \xymatrix{ H(G,\varphi) \ar[r]^{\varphi^s}\ar[d]^{\Psi} & H(G,\varphi)\ar[d]^\Psi\\
H(G',\varphi') \ar[r]^{\varphi'} & H(G',\varphi').}\]
\en
\end{lemma}

\begin{proof}
By the exact sequence \eqref{equ:mv} we have
$A'\cap B'=\Psi(I)$ where
\[ I=
\left\{ b\in B :\, \ba{l} \text{there exist $a\in A$ and $a_i\in A$, $i=0,\dots,s-2$, with}\\
                                        b\shtimes 0-a\shtimes (s-1)=\sum\limits_{i=0}^{s-2} a_i\shtimes i-\varphi(a_i)\shtimes (i+1)\in \bigoplus\limits_{i=0}^{s-1} G_i\ea \right\}. \]
\begin{claim}
$I=H(G,\varphi)$.
\end{claim}

\begin{proof}[Proof of the claim]
Let $b\in H(G,\varphi)$. By Lemma~\ref{lem:altdefh} we know that $\varphi^i(b)$ is defined for all $i$.
It is now straightforward to check that $a_i=\varphi^i(b)$, $i=0,\dots,s-2$ and $a=\varphi^{s-1}(b)$ (all of which lie in $A$, since $\varphi^i(b)$ is defined for any $i$) satisfy
\[ b\shtimes 0-a\shtimes (s-1)=\sum_{i=0}^{s-2} a_i\shtimes i-\varphi(a_i)\shtimes (i+1)\in \bigoplus\limits_{i=0}^{s-1} G_i,\]
that is, $b\in I$. On the other hand assume we have $b\in I$.
Let $a\in A$  and $a_i\in A$, $i=0,\dots,s-2$ as in the definition of $I$.
We deduce immediately that $b=a_0$, $a_{i+1}=\varphi(a_{i})$ for $i=0,\dots,s-2$ and $a=\varphi(a_{s-2})$. In particular
$\varphi^{i}(b)$ exists for $i=1,\dots,s-1$.
But by \eqref{equ:props}  this implies that $b\in H(G,\varphi)$.
\end{proof}

Put differently, the above claim shows that $\Psi$ gives rise to an isomorphism $H:=H(G,\varphi)\to A'\cap B'$.

\begin{claim}
The following diagram commutes:
\[ \xymatrix{ H\ar[d]^{\Psi}\ar[r]^{\varphi^s}&B\ar[d]^\Psi \\ A'\cap B'\ar[r]^{\varphi'}&B'.}\]
\end{claim}

\begin{proof}[Proof of the claim]
Let $b\in H$. The above discussion shows
that $\Psi(b)=b\times 0\in B'$ and $\varphi^{s-1}(b)\in A\times (s-1)=A'$ represent the same element in $G'$.
We now have
\[ \Psi(\varphi^s(b))=\varphi^s(b)\times 0=\varphi'(\varphi^{s-1}(b)\times (s-1))=\varphi'(\Psi(b)) \in G'.\]
This shows that the diagram commutes as claimed.
\end{proof}

Now note that  $\varphi^s$ defines an automorphism of $H$, hence $\varphi'$ defines an automorphism of $A'\cap B'$.
This shows that $A'\cap B'=H(G',\varphi')$. This concludes the proof of (1).
Statement (2) is now an immediate consequence from the above claims.
\end{proof}

We finally recall the following fact (cf.~\cite[Lemma~1.5]{Gr57}):

\begin{lemma}\label{lem:Gruenberg}
Let $G$ be a group and $N$ be a normal subgroup of $G$. If $G/N$ is a $p$-group and $N$ is residually $p$, then $G$ is residually $p$.
\end{lemma}

We are now in a position to complete the proof of Theorem~\ref{thm:respabelian} (and hence, of Theorem \ref{thm:resp}).
Suppose the order of $\varphi$ restricted to $H(G,\varphi)$ is a power of $p$. We continue to use the notations introduced above.
Pick  $k$ such that $p^k>r$ and set $s:=p^k$. Recall that $\phi$ denotes the morphism
$$G^*=\ll G,t \, |\, t^{-1}At=\varphi(A)\rr \to \Z/s\Z$$ with $t\mapsto 1$ and $g\mapsto 0$ for $g\in G$.
By Lemma~\ref{lem:Gruenberg} it suffices to show that $\G=\ker\phi$ is residually $p$.

By Lemma \ref{lem:abprime} we have
$H(G',\varphi')=A'\cap B'$.
Furthermore recall that we assumed that  the order of $\varphi|_{H(G,\varphi)}$ is a power of $p$. We can therefore appeal to  Lemma~\ref{lem:abprime} to see that  the order of the automorphism of $A'\cap B'$ induced by $\varphi'$ is a power of $p$.
Hence by Proposition~\ref{prop:hnnresp} and Lemma \ref{lem:resp} there exists a $p$-group  $X$ which contains $G'$ as a subgroup,
and an automorphism $\a\colon X\to X$ which extends the isomorphism $\varphi'\colon A'\to B'$ between subgroups of $G'$, and such that the order of $\a$ equals $p^l$ for some $l$.
We can therefore form the semidirect product $\Z/p^l\Z\ltimes X$ where $1\in \Z/p^l\Z$ acts on $X$ on the right via $\a$.
Now consider the composition $\psi$ of the morphisms
\[ \ba{rcl} \Gamma=\ll K,t \, |\, t^{-1}A_{s-1}t=B_{0}\rr \to \ll G',t \, |\, t^{-1}A't=B'\rr&\to& \Z/p^l\Z\ltimes X\\
t&\mapsto & (1,0)\\
g\in G'&\mapsto &(0,g).\ea \]
By Lemma~\ref{lem:giinj} the natural morphism $G_i\to G'$ is injective  for any $i\in \Z/p^k\Z$; hence the restriction of $\psi$ to $G_i$ is injective for any  $i\in \Z/p^k\Z$.

We are now finally in a position to prove that
 $\Gamma$ is residually $p$.
Recall that $\Gamma$ is formed as an HNN extension of the group
 $$K=G_0*_{A_0=B_1} G_1*_{A_1=B_2}G_2 \cdots \ast_{A_{s-2}=B_{s-1}} G_{s-1}.$$
One may think of $\Gamma$ as the fundamental group of a graph of groups with vertex groups $G_i$ ($i\in\Z/s\Z$) and edge groups $A_i=B_{i+1}$ ($i\in\Z/s\Z$).
 By the above the map $\psi\colon\Gamma\to G'$ is injective on any vertex group.
 It now follows from \cite[II.2.6, Lemma~8]{Se80} that $\ker(\psi)$ is a free group. Since free groups are residually $p$ and since $\ker(\psi)\leq \Gamma$
 is of $p$-power index it follows  from Lemma \ref{lem:Gruenberg} that $\Gamma$ is residually $p$.
 This concludes the proof of Theorem~\ref{thm:respabelian}. \qed

\end{document}